\documentclass{article}
\usepackage[utf8]{inputenc}
\usepackage{style}
\usepackage[margin=1in]{geometry}
\usepackage[english]{babel}
\usepackage{verbatim}
\usepackage{todonotes}

\usepackage[cal=cm]{mathalfa}

\DeclareMathOperator{\A}{\mathcal{A}}

\DeclareMathOperator{\B}{\mathcal{B}}

\usepackage{subfiles} % Best loaded last in the preamble

\title{Specialization of Integral Closure of Ideals by General Elements}
\author{Lindsey Hill and Rachel Lynn}
\date{ }

\begin{document}

\maketitle

\begin{abstract}

In this paper, we prove a result similar to results of Itoh and Hong-Ulrich, proving that integral closure of an ideal is compatible with specialization by a general element of that ideal for ideals of height at least two in a large class of rings. 
Moreover, we show integral closure of sufficiently large powers of the ideal is compatible with specialization by a general element of the original ideal.
In a polynomial ring over an infinite field, we give a class of squarefree monomial ideals for which the integral closure is compatible with specialization by a general linear form. 

\end{abstract}

\section{Introduction}

\indent In this paper we prove that for ideals of height at least two in a large class of rings, the integral closure of the ideal is compatible with specialization modulo general elements of the ideal. 

Let $R$ be a ring, and let $I$ be an $R$-ideal. 
The Rees algebra of $I$, denoted $R[It]$, is the ring 
\[R[It] = \bigoplus_{n = 0}^{\infty} I^nt^n \subset R[t]\]
where $t$ is an indeterminate over the ring $R$. 
The integral closure of the Rees algebra in $R[t]$, which we will denote as $\overline{R[It]}$, is a graded algebra and has the form \[\overline{R[It]} =  \bigoplus_{n = 0}^{\infty} I_n t^n,\] 
where each $I_n$ is an $R$-ideal. $I_1$ is called the integral closure of $I$, and is denoted $\overline{I}$.
Equivalently, we can define the integral closure of an ideal to be the elements of $R$ which satisfy an equation of integrality over $I$: $\overline{I} = \{x\in R \,| \, x^n+a_1 x^{n-1} + \cdots + a_n = 0\text{, with } a_i \in I^i\}$. 
In fact, the integral closure of $I^n$ is the degree $n$ component of $\overline{R[It]}$.  

For more information about integral closures, including a proof of the form of $\overline{R[It]}$, we refer the reader to \cite{Swanson-Huneke}.

In this paper, we are interested in determining whether the integral closure of an ideal is compatible with taking the quotient by an element of the ring; we call this quotient by an element \textit{specialization}. 
In our main result, \cref{thm}, we prove that for most algebras over a field of characteristic zero, specialization by a general element of the ideal (which we define in Section 2) is compatible with integral closure for ideals of height at least two; that is, $\overline{I}/(x) = \overline{I/(x)}$ for a general element $x$ of $I$.

If we choose a specific element to specialize by, it is certainly not the case that being integrally closed is preserved. 
Consider, for example, the two-dimensional polynomial ring $R=k[x,y]$ over a field $k$, and the integrally closed ideal $I = \overline{I} = \left(x^2,xy,y^2\right)$.
Going modulo the generator $x^2$, we see that $x$ is integral over $I/(x^2)$ since it satisfies the equation of integral dependence $z^2=0$ in $R/(x^2)$. 
Since $x \not\in I/(x^2)$, we see that $I/(x^2)$ is not integrally closed.  
However, in \cite{Itoh-Coefficients} Itoh proved specialization is compatible with integral closures for parameter ideals in analytically unramified local Cohen-Macaulay rings. 
Hong and Ulrich then generalized this result to a larger class of ideals in a larger class of rings in \cite{Hong-Ulrich}, proving that being integrally closed is preserved under specialization by generic elements of $I$ if $R$ is a Noetherian, locally equidimensional, universally catenary ring such that $R_{\operatorname{red}}$ is locally analytically unramified and the height of $I$ is at least 2. 
These results have many applications, including a simple proof of the theorem (proved independently by Itoh in \cite{Itoh} and Huneke in \cite{Huneke}) that for an ideal $I$ generated by a regular sequence $x_{1},\dots, x_{r}$ in a Noetherian local ring, $I^{n} \cap \overline{I^{n+1}} = I^{n} \overline{I}$; 
Huneke's proof of this theorem partially inspired the creation of tight closure. 
Hong and Ulrich's result can also be used to show the integral closure of a module $M$ is compatible with specialization by a generic element of the module through the use of the generic Bourbaki ideal. 

The results of Hong-Ulrich and Itoh prove results for specialization by a generic element, which requires extending the original ring $R$ to a polynomial ring over $R$, or a localization thereof, in order to obtain the generic element. 
This is often sufficient when one wishes to use the result in a proofs using induction on the height of the ideal when considering properties preserved under faithfully flat extensions. 
However, there are cases where we may not want to change the ring, such as when considering the core of the ideal or when it is important that the residue field is perfect or algebraically closed, since these properties may not pass to extensions of $R$.
To preserve the original ring we need to instead consider \textit{general} elements of the ideal, which come from the original ring.
An important breakthrough for specialization by general elements came from Flenner's local version of Bertini's theorem in \cite{Flenner}, which tells us how Serre's conditions pass to the specialization of a ring.

In Section 2, we discuss results from the literature, such as Flenner's theorem, which are important to the proof of our main theorem. 
In Section 3, we provide the proof of our main theorem, which proves that integral closure and specialization by a general element of the ideal are compatible for a large class of ideals and rings. We also prove that integral closure and specialization by a general element of the ideal are compatible for sufficiently large powers of the ideal. 
In Section 4, we consider the question of specialization by a general element of the unique maximal ideal rather than a general element of the ideal $I$ and prove that specialization by a general linear form is compatible with integral closure for a class of squarefree monomial ideals.

 \label{intro}

\section{Preliminaries}

In this section, we introduce results from the literature that are helpful in proving our main theorem. First, we provide a few definitions.

\begin{definition}\label{def-extended-rees-algebra}
Let $R$ be a ring, and let $I$ be an $R$-ideal. The \textbf{extended Rees algebra} of $I$, denoted $R[It, t^{-1}]$, is the $R$-subalgebra of the Laurent polynomial ring $R[t, t^{-1}]$
\[
R[It, t^{-1}] = \bigoplus_{n = - \infty}^{\infty} I^n t^n \subset R[t, t^{-1}],
\] 
where $t$ is an indeterminate of the ring $R$.
\end{definition}

The integral closure of the extended Rees algebra in the Laurent polynomial ring, which we will denote as $\overline{R[It, t^{-1}]}$, is a graded algebra and has the form
\[
\overline{R[It, t^{-1}]} = \bigoplus_{n = - \infty}^{\infty} \overline{I^n} t^n,
\]
where $\overline{I^n}$ is defined to be equal to $R$ for $n \leq 0$. Once again, we refer the reader to \cite{Swanson-Huneke} for more information.

\begin{definition}\label{def-generic-element}
Let $R$ be a ring, let $I = (a_1,\dots, a_n)$ be an ideal generated by $n$ elements, and let $S = R[X_1,\dots,X_n]$ be a polynomial ring over $R$ in $n$ variables. If $(R,m)$ is a local ring, we may alternatively let $S = R[X_1,\dots, X_n]_{mR[X_1,\dots,X_n]}$ be the localization of a polynomial ring. An element $a = \sum_{i=1}^n X_i a_i $ is called a \textbf{generic element} of $I$. 
\end{definition} 

Notice that a generic element of an ideal $I$ lies in the ideal $IS$, in the larger ring $S$, rather than in the original ring $R$; thus the results of Itoh and Hong-Ulrich discussed in \cref{intro} may not apply when the original ring must be preserved.

We now define a general element of an ideal $I$ in an algebra over an infinite field. A general element is similar to a generic element, but sits inside of the original ring.

\begin{definition}\label{def-general-element}
Let $k$ be an infinite field, let $R$ be a $k$-algebra, and let $I = (a_1,\dots, a_n)$. For $\alpha = (\alpha_1,\dots, \alpha_n) \in k^n$, let $a_\alpha = \sum\limits_{i=1}^n \alpha_i a_i$. Then a property $P$ holds for a \textbf{general element} of $I$ if there is a dense open subset $U$ in affine $n$-space $k^n$ such that $P$ holds for $a_\alpha$ whenever $\alpha\in U$.
\end{definition}

\begin{remark}
We can extend the definition of a general element to an ideal $I = (a_1, \dots, a_n)$ of a local ring $(R,m)$ which has infinite residue field $k$. For $\alpha = (\alpha_1,\dots, \alpha_n) \in R^n$, we say a property $P$ holds for a \textbf{general element} of $I$ if there is a dense open subset $U$ of affine $n$-space $k^n$ such that the property $P$ holds for $a_\alpha = \sum_{i=1}^n \alpha_i a_i$ whenever the image $\overline{\alpha}$ of $\alpha = (\alpha_1,\dots,\alpha_n)\in R^n$ in $k^n$ belongs to $U$.
\end{remark}

\begin{remark}\label{rmk-finite-intersection}
If $k$ is an infinite field, the intersection of finitely many dense open subsets of $k^n$ is again a dense open subset of $k^n$. It follows that if finitely many properties $P_1,\dots, P_s$ hold for a general element of $I$, then the intersection of all the properties, $P_1\wedge \cdots \wedge P_s$, also holds for a general element of $I$.
\end{remark}

Next we note a lemma of Hong and Ulrich which they used to prove their result on the compatibility of specialization by generic elements with integral closure. The statement follows from {\cite[Lemma~1.1(c)]{Hong-Ulrich}}, although it is not the stated result.

\begin{lemma}[{\cite[Lemma~1.1(c)]{Hong-Ulrich}}]\label{1.1}\label{c}
Let $R$ be a Noetherian, equidimensional, universally catenary local ring of dimension $d$ such that $R_{\operatorname{red}} = R/\sqrt{0}$ is analytically unramified. Let $I = (a_1,\dots, a_n)$ be a proper $R$-ideal with $\hgt{I} > 0$ and write $\A = R[It, t^{-1}]$ for the extended Rees ring of $I$.
Then  $\grade{It}{\intcl{\A}/t^{-1}\intcl{\A}} > 0.$
\end{lemma}

The next lemma follows directly from \cref{c} above:
\begin{lemma}\label{depth}
Using the notation of \cref{1.1}, let $x$ be a general element of $I$. Then we may assume $t^{-1}, xt$ is a regular sequence on $\overline{\A}$. In particular, $xt$ is regular on $\overline{\A}/t^{-1}\overline{\A}$.
\end{lemma}

In \cref{flenner47} we cite a theorem of Flenner. In his paper, Flenner proves there is a local version of Bertini's second classical theorem, positively answering a question posed by Grothendieck. This result is vital to the proof of our main theorem, as it allows us to pass Serre's conditions to specializations of $\overline{R[It]}$.

\begin{theorem}[{\cite[Corollary~4.7]{Flenner}}]\label{flenner47}
Let $(S, m_S)$ be a local excellent $k$-algebra over a field $k$ of characteristic 0, let $I = (x_1,\dots, x_n) \subset m_S$, let $V(I) = \{p\in \operatorname{Spec}{S} \st I\subseteq p \}$, and let $D(I)=\operatorname{Spec}(S)\setminus V(I)$. Assume that for all $p\in D(I)$, the ring $S_p$ satisfies Serre's condition $(S_r)$ or $(R_s)$. For $\alpha \in k^n$, let $x_\alpha = \sum_{i=1}^n \alpha_i x_i$. Then for general $\alpha \in k^n$ and $p\in D(I)\cap V(x_\alpha)$ the ring $(S/x_\alpha S)_p$ also satisfies the condition $(S_r)$ or $(R_s)$.
\end{theorem}

In particular, normality passes to the specialization of a ring in the following way:

\begin{cor}[{\cite[Corollary~4.8]{Flenner}}]\label{flenner}
Let $S$ be a local excellent $k$-algebra over the field $k$ of characteristic 0 and let $I = (x_1,\dots, x_n) \subset m_S$. Let $\operatorname{Nor}(S)=\{p \in \operatorname{Spec}(S) \mid S_{p} $ is normal$\}$, For general $\alpha \in k^n$, let $x_{\alpha}= \sum_{i=1}^n \alpha_i x_i$, as in \cref{flenner47}. Then
\[\operatorname{Nor}(S)\cap V(x_\alpha)\cap D(I) \subset \operatorname{Nor}(S/x_\alpha S).\]
\end{cor}

A final ingredient in the proof of \cref{thm} is the vanishing of local cohomology of the integral closure of the extended Rees ring of $I$.

\begin{theorem}[{\cite[Theorem~1.2]{Hong-Ulrich}}]\label{1.2}
Let $R$ be a Noetherian, locally equidimensional, universally catenary ring such that $R_{\operatorname{red}}$ is locally analytically unramified. Let $I$ be a proper $R$-ideal with $\hgt{I}>0$, $\A=R[It,t^{-1}]$ the extended Rees ring of $I$, and $\overline{\A}$ the integral closure of $\A$ in $R[t,t^{-1}]$. Let $J$ be an $\A$-ideal of height at least $3$ generated by $t^{-1}$ and homogeneous elements of positive degree. Then $[H_{J}^{2}(\overline{\A})]_{n}=0$ for all $n \leq 0$, where $H_{J}^i(\phantom{A})$ denotes the $i^{th}$ local cohomology with respect to the ideal $J$, and $[\phantom{m}]_{n}$ denotes the degree $n$ component.
\end{theorem}

\section{Specialization by a General Element of $I$}

In this section we prove that integral closure of an ideal is compatible with specialization by a general element of the ideal. First, we state and prove a few technical lemmas necessary to reduce to the case where $R$ is a normal reduced ring. \cref{red} allows us to restrict our focus to reduced rings, \cref{norm} allows us to restrict our focus further to normal rings, and \cref{hgtI} shows the height of an ideal is preserved under these reductions.
%-------------------------------------------------------------------------------------------

\begin{lemma}\label{red}
Let $R$ be an algebra over a field $k$ with $|k| = \infty$, let $R_{\operatorname{red}}=R/\sqrt{0}$, and let $J$ be an $R$-ideal. Let $x$ be a general element of an $R$-ideal $I$, where $I$ may or may not be equal to $J$. If specialization by $x$ is compatible with integral closure for the image of $J$ in $R_{\operatorname{red}}$, then specialization and integral closure are also compatible for the ideal $J$. Stated symbolically, if  $\overline{JR_{\operatorname{red}}}\left(R_{\operatorname{red}}/(x)R_{\operatorname{red}}\right)=\overline{J\left(R_{\operatorname{red}}/(x)R_{\operatorname{red}}\right)}$, then $\overline{J}\left(R/(x) \right)=\overline{J \left(R/(x) \right)}$.
\end{lemma}

\begin{proof}
The containment $\overline{J}\left(R/(x) \right) \subseteq \overline{J \left(R/(x) \right)}$ is immediate by persistence of the integral closure.
It remains to show the reverse containment.

Let $\psi$ denote the natural map from $R/(x)$ to $R_{\operatorname{red}}/(x)R_{\operatorname{red}}$. 
Then by assumption 
\[
\psi^{-1}\left(\overline{JR_{\operatorname{red}}}\left(R_{\operatorname{red}}/(x)R_{\operatorname{red}}\right)\right) =\psi^{-1}\left(\overline{J\left(R_{\operatorname{red}}/(x)R_{\operatorname{red}}\right)}\right).\] 
By \cite[Proposition~1.1.5]{Swanson-Huneke}, $\overline{JR_{\operatorname{red}}}=\overline{J}R_{\operatorname{red}}$, which implies that \begin{align*}
\psi^{-1}\left(\overline{JR_{\operatorname{red}}}\left(R_{\operatorname{red}}/(x)R_{\operatorname{red}}\right)\right)
&=\psi^{-1}(\overline{J}\left(R_{\operatorname{red}}/(x)R_{\operatorname{red}}\right)) \\ 
&= \{a+(x) \mid a+(x)+\sqrt{0} \in \overline{J}+(x)+\sqrt{0})\}\\
&= \{a+(x) \mid a+(x) \in \overline{J}+(x))\} \text{ because $\sqrt{0}\subset \overline{J}$}\\
&=\overline{J}\left(R/(x)\right).
\end{align*}
Hence it suffices to show that $\overline{J \left(R/(x)\right)} \subset \psi^{-1}(\overline{J \left(R_{\operatorname{red}}/(x)R_{\operatorname{red}}\right)})$. 
By persistence and taking preimages, we see that \[\overline{J\left(R/(x)\right)}\subset \psi^{-1} \left( \psi\left( \overline{J\left(R/(x)\right)} \right) \right)\subset \psi^{-1}\left(\overline{J \left(R_{\operatorname{red}}/(x)R_{\operatorname{red}}\right)}\right).\] 
Therefore, we conclude that $\overline{J \left(R/(x)\right)}=\overline{J}\left(R/(x)\right).$
\end{proof}

\begin{lemma}\label{norm}
Let $R$ be an algebra over a field $k$ with $|k| = \infty$, let $m\in \mspec{\overline{R}}$, let $I$ be an $R$-ideal, and let $x$ be a general element of $I$. If specialization by $x$ is compatible with integral closure for the image of $I$ in $\overline{R}_{m}$ for every maximal ideal $m$, then specialization and integral closure are also compatible for the ideal $I$. Stated symbolically, if $\overline{I\overline{R}_{m}}\left(\overline{R}_{m}/(x)\overline{R}_{m}\right) = \overline{I\left(\overline{R}_{m}/(x)\overline{R}_{m}\right)}$ for all $m\in \mspec{\overline{R}}$, then $\overline{I}\left(R/(x)\right)=\overline{I\left(R/(x)\right)}$. 
\end{lemma}

\begin{proof}
We first show that $\overline{I\overline{R}}\left(\overline{R}/(x)\overline{R}\right)$ $=\overline{I\left(\overline{R}/(x)\overline{R}\right)}$. 
It suffices to check equality locally at maximal ideals $m \in \mspec{\overline{R}/(x)\overline{R}}$. Identifying $\mspec{\overline{R}/(x)\overline{R}}$ with a subset of $\mspec{\overline{R}}$, by compatibility of integral closure and localization and by assumption, we have the following equalities: 
\[ \left(\overline{I \left(\overline{R}/(x)\overline{R}\right)}\right)_{m}
=\overline{I\left(\overline{R}_{m}/(x)\overline{R}_{m}\right)}
=\overline{I\overline{R}_{m}}\left(\overline{R}_{m}/(x)\overline{R}_{m}\right)
=\left(\overline{I\overline{R}}\left(\overline{R}/(x)\overline{R}\right)\right)_{m}.\]
This proves the desired equality.

By persistence of integral closure, $\overline{I}\left(R/(x)\right) \subset \overline{I\left(R/(x)\right)}$. It remains to show the reverse containment. Let $\psi$ denote the natural map from $R/(x)$ to $\overline{R}/(x) \overline{R}$. Then by above, \[\psi^{-1}\left(\overline{I\overline{R}}\left(\overline{R}/(x)\overline{R}\right)\right) 
= \psi^{-1}\left(\overline{I\left(\overline{R}/(x)\overline{R}\right)}\right).\] By \cite[Proposition~1.6.1]{Swanson-Huneke}, $\overline{I\overline{R}} \cap R=\overline{I}$, thus $\psi^{-1}\left(\overline{I\overline{R}}\left(\overline{R}/(x)\overline{R}\right)\right)
=\overline{I}\left(R/(x)\right)$. Therefore, it suffices to show that $\overline{I\left(R/(x)\right)} \subset
\psi^{-1}\left(\overline{I\left(\overline{R}/(x)\overline{R}\right)} \right)$. By persistence and taking preimages, we see that 
\[ \overline{I\left(R/(x)\right)} \subset
\psi^{-1}\left(\psi\left(\overline{I\left(R/(x)\right)}\right)\right) \subset \psi^{-1}\left(\overline{I\left(\overline{R}/(x) \overline{R}\right)}\right).\]
Therefore, we conclude that $\overline{I}\left(R/(x)\right)=\overline{I\left(R/(x)\right)}$.

\end{proof}

\begin{lemma} \label{hgtI} 
Let $(R, m)$ be a local equidimensional excellent ring. Let $I$ be an $R$-ideal. Then $\hgt{I\overline{R_{\operatorname{red}}}} = \hgt{I}$.
\end{lemma}

\begin{proof}
Since $\hgt{I} = \hgt{IR_{\operatorname{red}}}$, we may pass from $R$ to $R_{\operatorname{red}}$ to assume $R$ is a reduced local equidimensional excellent ring. It remains to show that $\hgt{I\overline{R}} = \hgt{I}$. 

We first note that $\overline{R}$ is finitely generated over $R$ since $R$ is a reduced excellent ring.
Since $R$ is excellent, it is by definition universally catenary. Therefore, $\intcl{R}$ is catenary. 

We claim that there is a one-to-one correspondence of minimal primes of $R$ and minimal primes of $\overline{R}$. Let $S$ denote the set of nonzerodivisors on $R$, and let $W$ denote the set of nonzerodivisors on $\overline{R}$. 
Since $R \subseteq \overline{R}$ and $S \subseteq W$, 
\begin{equation*}
   \Quot{R}= S^{-1}R \subseteq S^{-1}\overline{R} \subseteq W^{-1}\overline{R}= \Quot{\overline{R}}.
\end{equation*}

Next, we see that every element of $W$ is a unit in $\Quot{R}$. Let $w \in W$. Then since $W \subseteq \overline{R} \subseteq \Quot{R}$, $w=u/v$ for some $u \in R$ and $v\in S$. 
Moreover, since $u/v$ is a nonzerodivisor on $\overline{R}$, $u$ is a nonzerodivisor on $R$. Hence $w=u/v$ is a unit in $\Quot{R}$. 
Then since $\overline{R} \subseteq \Quot{R}$ and $W \subseteq (\Quot{R})^{*}$, the units of $\Quot{R}$, we conclude that 
\begin{equation*}
\Quot{\overline{R}} = W^{-1}\overline{R} \subseteq W^{-1}\Quot{R} = \Quot{R}.
\end{equation*}

This shows that $\Quot{R}=\Quot{\overline{R}}$; i.e., $R \rightarrow \overline{R}$ is a birational extension. 

Note that for any Noetherian ring $T$, since the total ring of quotients $\Quot{T}$ is the localization of $T$ with respect to the complement of the union of the associated primes of $T$, the minimal primes of $T$ correspond to the minimal primes of $\Quot{T}$. Since the minimal primes of $R$ and $\overline{R}$ are both in one-to-one correspondence with the minimal primes of $\Quot{R}=\Quot{\overline{R}}$, we conclude that $\Min{R}$ is in one-to-one correspondence with $\Min{\overline{R}}$. Therefore, every minimal prime of $\overline{R}$ contracts to a minimal prime of $R$. This proves our claim.

Now let $m\in\mspec{\overline{R}}$. 
Let $q \in \Min{\overline{R}}$ be contained in $m$. 
Notice $m\cap R$ must be the unique maximal ideal of $R$ since $R \rightarrow \overline{R}$ is an integral extension, and $q\cap R$ must be a minimal prime of $R$.
Since $R$ is equidimensional and local, $\dim{\left(R/(q\cap R)\right)_{m\cap R}} = \dim{R}$. 
Notice that the dimension formula applies because $R/(q\cap R) \subset \intcl{R}/q$ is an extension of domains, $R/(q\cap R)$ is a universally catenary Noetherian ring, and $\intcl{R}/q$ is a finitely generated algebra over $R/(q\cap R)$.
By the dimension formula, 
\begin{align*}
    \dim{\left(\overline{R}/q \right)_{m}} &= \dim{\left(R/(q\cap R)\right)_{m\cap R}} + \operatorname{trdeg}_{R/(q\cap R)} \overline{R}/q  - \operatorname{trdeg}_{\kappa\left((m\cap R)/(q\cap R)\right)} \kappa(m/q) \\
    & = \dim{R} + 0 - 0 = \dim{\overline{R}},
\end{align*}
thus $\overline{R}$ is locally equidimensional of the same dimension at every maximal ideal $m$. 

Since $\overline{R}$ is catenary and locally equidimensional of the same dimension at every maximal ideal $m$, one sees that for any prime $p\in\Spec{\overline{R}}$, $\dim{\overline{R}}/p + \hgt{p} = \dim{\overline{R}}$, and hence the same holds for any ideal of $\overline{R}$. 
Thus 
\begin{align*}
    \dim{R/\overline{I}} + \hgt{\overline{I}} 
    &=\dim{R}\\
    &=\dim{\overline{R}}\\
    &= \dim{\overline{R}/\overline{I\overline{R}}} + \hgt{\overline{I\overline{R}}}\\
    &= \dim{R/{\overline{I}}} + \hgt{\overline{I\overline{R}}} \text{\hspace{.25in}since $\overline{I\overline{R}}\cap R = \overline{I}$ by \cite[Proposition~1.6.1]{Swanson-Huneke}.}
\end{align*}

From the above equalities we see that $\hgt{\overline{I}} = \hgt{\overline{I\overline{R}}}$.
We conclude that $\hgt{I} = \hgt{\overline{I}} = \hgt{\overline{I\overline{R}}} = \hgt{I\overline{R}}$ since an ideal and its integral closure have the same minimal primes.
\end{proof}

We are now ready to prove our main theorem.

\begin{theorem}\label{thm}
Let $(R, m)$ be a local equidimensional excellent $k$-algebra, where $k$ is a field of characteristic 0. 
Let $I =(a_1,...,a_n)$ be an $R$-ideal such that $\hgt I \geq 2$, and let $x$ be a general element of $I$. 
Then $\overline{I}/(x)= \overline{I/(x)}$.
\end{theorem}

\begin{proof}
By \cref{red}, we may pass from $R$ to $R_{\operatorname{red}}$ to assume $R$ is a reduced local equidimensional excellent $k$-algebra and the height of $I$ does not change. 
Since $R$ is a reduced excellent local ring, $\intcl{R}$ is finitely generated over $R$.
This shows that $\intcl{R}_m$ is excellent and that $\intcl{R}$ is semilocal.
Hence $x$ is still a general element of $I\intcl{R}_m$ for each of the finitely many maximal ideals $m$ of $\intcl{R}$. 
By \cref{norm} we may pass from $R$ to $\overline{R}_{m}$ for any $m \in \mspec{\overline{R}}$ to assume in addition that $R$ is a local normal ring, hence also a domain, and by \cref{hgtI} we may still assume that $I$ has height at least $2$. 

Let $\A$ and $\B$ denote the extended Rees algebras of $I$ and $I/(x)$, respectively.
Let $\overline{\A}$ and $\overline{\B}$ denote the integral closures of $\A$ and $\B$ in $R[t,t^{-1}]$ and $\left(R/(x)\right)[t,t^{-1}]$ respectively.
Notice that because $R$ is excellent, $\intcl{\A}$ is finitely generated over $\A$ and is thus Noetherian. 
Define $J$ to be the $\A$-ideal $\left(It, t^{-1}\right)$. 
We consider the natural map $\varphi:\overline{\A}/xt\overline{\A} \rightarrow \overline{\B}$ induced by $R[t, t^{-1}] \rightarrow (R/(x))[t,t^{-1}]$. 
It suffices to show that $\coker{\phi}$ in degree 1 is 0. 
We first show that $\varphi_{p}$ is an isomorphism for $p \in \operatorname{Spec}({\overline{\A}}) \setminus V(J\overline{\A})$. 

First suppose $t^{-1}\not\in p$. 
Then one sees that
\[ \left(\overline{\A}/xt\overline{\A} \right)_{t^{-1}} \cong \frac{\overline{\A}_{t^{-1}}}{xt\overline{\A}_{t^{-1}}} \cong \frac{R[t,t^{-1}]}{xtR[t,t^{-1}]} = R[t, t^{-1}]/xR[t,t^{-1}] \cong \frac{R}{(x)}[t,t^{-1}] \cong \overline{\B}_{t^{-1}}.\] 
Therefore by further localization, $\left(\overline{\A}/xt\overline{\A} \right)_{p} \cong \overline{\B}_{p}$ if $t^{-1} \not\in p$.

Now let $p \in \operatorname{Spec}({\overline{\A}}) \setminus V(It\overline{\A})$. 
We first show $\phi_p$ is injective. Notice \[(xt\overline{\A})_p\subset (xtR[t,t^{-1}]\cap \overline{\A})_p\] and \[ (xtR[t,t^{-1}]\cap \overline{\A})_p/xt\intcl{\A}_p = \ker\phi_p,\] so it is enough to show the two are equal locally at associated primes of $xt\overline{\A}_p$. 
Additionally, $\overline{\A}_{p}$ is normal and excellent since $R$ is normal and excellent. 
Thus by \cref{flenner}, $(\overline{\A}/xt\overline{\A})_p$ is normal, hence a domain. 
Thus $xt\overline{\A}_p$ is a prime ideal, and it suffices to show equality locally at $q = xt\intcl{\A}_p$. 
Notice $\hgt{q} \leq 1$. 
Since $t^{-1}, xt$ is an $\overline{\A}$-regular sequence by \cref{depth}, we may assume $t^{-1} \not \in q$. 
Thus the desired equality holds locally at $q$ by the previous case. 
Thus $\phi_p$ is injective for all $p\in \operatorname{Spec}({\overline{\A}}) \setminus V(It\overline{\A})$.

We now show that for the same primes $p$, $\phi_{p}$ is also a surjection. 
Since $\A$ surjects onto $\B$ we have that the ring extension $\im(\A)_{p\cap \A} = \B_{p\cap \A} \subset\overline{\B}_{p\cap \A}$ is an integral extension. 
Hence $\im(\overline{\A})_{p\cap \A} \subset \overline{\B}_{p\cap \A}$ is an integral extension, 
and thus $\im(\overline{\A})_p = \im \left(\overline{\A}/xt\overline{\A}\right)_p \subset \overline{\B}_p$ is also an integral extension. 

Notice $\im(\left(\intcl{\A}/xt\intcl{\A}\right)_p) \subset \intcl{\B}_p$, $t^{-1}$ is a nonzerodivisor on both rings, and after making $t^{-1}$ invertible both rings are equal. 
Hence they have the same total ring of quotients.

Recall that $\left(\overline{\A}/xt\overline{\A}\right)_p$ is normal and $\phi_p$ is injective. 
Hence $\im(\intcl{\A}/xt\intcl{\A})_p$ is also normal, so $\phi_p$ must be a surjection. 
Therefore $\phi_p$ is an isomorphism, proving our claim.

% PART 2----------------------------------------------------------------------------------
Let $K$ denote the kernel of $\varphi$ and $C$ the cokernel of $\varphi$. 
In order to show that $[C]_{1}=0$, we will identify $C$ with a submodule of the local cohomology module $H_J^2(\overline{\A})$. 
By the isomorphisms above, for all $p\not\in V(J\overline{\A})$ we have $K_p = C_p = 0$.
Hence $H_J^0(K) = K$ and thus $H_J^i(K) = 0$ for all $i>0$. Similarly $H_J^0(C) = C$. 
We also have $H_J^0(\overline{\B}) = 0$ because $t^{-1} \in J$ is a $\overline{\B}$-regular element. 

From the long exact sequence of local cohomology induced by the exact sequence 
\[
\begin{tikzcd}
0\rar & K \rar & \overline{\A}/xt\overline{\A} \rar{\phi} & \operatorname{Im}(\phi) \rar & 0
\end{tikzcd}
\]
we obtain $H_J^i(\overline{\A}/xt\overline{\A}) \cong H_J^i(\operatorname{Im}(\phi))$ for all $i\geq 1$ since $H_J^i(K)$ vanishes for $i\geq 1$.

From the long exact sequence of local cohomology induced by the exact sequence 
\[
\begin{tikzcd}
0\rar & \operatorname{Im}(\phi) \rar & \overline{\B} \rar & C \rar & 0
\end{tikzcd}
\]
we obtain the exact sequence
\[
\begin{tikzcd}
0 = H_J^0(\overline{\B}) \rar & H_J^0(C) \rar & H_J^1(\operatorname{Im}(\phi)) \cong H_J^1(\overline{\A}/xt\overline{\A}).
\end{tikzcd}
\]
Therefore $C = H_J^0(C) \hookrightarrow H_J^1(\overline{\A}/xt\overline{\A})$.

By \cref{depth} we have $\operatorname{depth}_J(\overline{\A}) \geq 2$. Thus $H_J^1(\overline{\A}) = 0$. 
Hence from the short exact sequence 
\[
\begin{tikzcd}
0\rar & xt\overline{\A} \rar & \overline{\A} \rar & \overline{\A}/xt\overline{\A} \rar & 0
\end{tikzcd}
\]
we obtain the exact sequence
\[
\begin{tikzcd}
0 = H_J^1(\overline{\A}) \rar & H_J^1(\overline{\A}/xt\overline{\A}) \rar & H_J^2(xt\overline{\A}).
\end{tikzcd}
\]
Therefore $C\hookrightarrow H_J^2(xt\overline{\A})$.

Since $x$ is a general element of $I$, we may assume $x$ is a nonzerodivisor on $R$, and hence that $xt$ is a nonzerodivisor on $\overline{\A}\subset R[t,t^{-1}]$. 
Thus we have an isomorphism of graded modules $xt\overline{\A} \cong \overline{\A}(-1)$, and hence $C \hookrightarrow H_J^2({\overline{\A}(-1)})$.
So $[C]_n \hookrightarrow [H_J^2({\A(-1)})]_n \cong [H_J^2({\overline{\A}})]_{n-1}$. 
By \cref{1.2}, we have $[H_J^2({\overline{\A}})]_{n-1} = 0$ for all $n \leq 1$ as long as $\hgt{J} \geq 3$. 
In particular, this will imply that $[C]_1 = 0$, that is, $\overline{I/(x)} = \overline{I} / (x)$.

It remains to show that $\hgt{J} \geq 3$.

Since $R$ is excellent and hence universally catenary, $\A$ is catenary. 
Since $R$ is equidimensional, the localization of $\A$ at the unique maximal homogeneous ideal $\mathfrak{m}$ is equidimensional of dimension $\dim{R} + 1$ by \cite[Theorem~5.1.4(3)]{Swanson-Huneke}. 

Since $J$ is a homogeneous $\A$-ideal, the minimal primes of $J$ are homogeneous. 
Hence $\hgt J= \hgt J_{\mathfrak{m}}$. 
Since $\A_{\mathfrak{m}}$ is equidimensional and catenary, we have that
\begin{equation*}
    \hgt{J} = \hgt{J_{\mathfrak{m}}}= \dim \A_{\mathfrak{m}}- \dim (\A/J)_{\mathfrak{m}}.
\end{equation*}
Since $J$ is homogeneous, $\dim (\A/J)_{\mathfrak{m}}$ is equal to $\dim (\A/J)$. 
Notice that $\A/J \cong R/I$ and therefore, 
\begin{align*}
    \hgt{J} &= \dim R + 1 - \dim A/J \\
    &= \dim R + 1 - \dim R/I \\
    &= \dim R + 1 -\dim R + \hgt{I} \\
    &= 1+\hgt{I} \\
    & \geq 3.
\end{align*}

Thus by \cref{1.2}, we have $[H_J^2({\overline{\A}})]_{n-1} = 0$ for all $n \leq 1$, and so $\overline{I/(x)} = \overline{I} / (x)$.

\end{proof}

We note that the assumption on the height of the ideal in \cref{thm} is necessary. We give a counterexample for height $1$ ideals. 

\begin{example}
Let $R=\mathbb{Q}[t]_{(t)}$, $I=\left( t^2 \right)$ and let $x$ denote a general element of $I$. Since $I$ is integrally closed and $(x)=I$, $\overline{I}/(x)$ is the zero ideal in $R/(x)$. However, $t+(x)$ is a nonzero element of $\overline{I/(x)}$, so $\overline{I}/(x)\subsetneq \overline{I/(x)}$. 
\end{example}

When $R$ is normal, we obtain a result similar to \cref{thm} for sufficiently large powers of the ideal.

\begin{remark}
In \cref{powers} we add the assumption that $R$ is normal. 
This is because we are unable to reduce to the normal case as in the proof of \cref{thm} since the proof of \cref{norm} does not work when $x$ is not an element of the ideal we are taking the integral closure of, in this case $I^s$. 
\end{remark}

\begin{prop}\label{powers}
Let $(R, m)$ be a local normal excellent $k$-algebra, where $k$ is a field of characteristic 0. 
Let $I =(a_1,...,a_n)$ be an $R$-ideal such that $\hgt I \geq 2$, and let $x$ be a general element of $I$. 
Then $\overline{I^{s}}+(x)/(x)= \overline{(I/(x))^{s}}$ for $s$ sufficiently large.
\end{prop}

\begin{proof}
By \cref{red}, we may reduce to the case that $R$ is in addition a reduced ring. 
Let $\mathcal{A}$ and $\mathcal{B}$ denote the extended Rees algebras of $I$ and $I/(x)$, respectively. 
Let $J$ denote the $\mathcal{A}$-ideal $(It, t^{-1})$. 
Consider the natural map $\varphi: \overline{\mathcal{A}}/xt\overline{\mathcal{A}} \rightarrow \overline{\mathcal{B}}$.
Denote the cokernel of $\varphi$ by $C$. As in the proof of \cref{thm}, $H_J^0(C)=C$. 

Next we show that we may assume $R/(x)$ is a reduced excellent ring. 
Since $R$ is excellent, $R/(x)$ is excellent. 
Since $R$ is reduced, by \cref{flenner47} $R/(x)$ satisfies $R_{0}$ locally at primes which do not contain $I$. 
But since $I$ has height at least $2$, $R/(x)$ satisfies $R_{0}$. 
Moreover, since $R$ satisfies $S_{2}$ and $x$ is a nonzerodivisor, $x$ can be chosen generally so that $R/(x)$ satisfies $S_{1}$.
Therefore, $R/(x)$ is reduced. Now, notice that $\mathcal{\B} \subset (R/(x))[t,t^{-1}]$ and therefore is also a reduced ring. 
Furthermore, $\mathcal{\B}$ is a finitely generated $R/(x)$-algebra and hence is excellent, so $\overline{\mathcal{B}}^{\text{Quot}(\mathcal{B})}$ is a finite $\mathcal{B}$-module.
Since $\mathcal{B}$ is a Noetherian ring, $\overline{\mathcal{B}}^{\text{Quot}(\mathcal{B})}$ is a Noetherian $\mathcal{B}$-module and therefore, $\overline{\mathcal{B}} \subset \overline{\mathcal{B}}^{\text{Quot}(\mathcal{B})}$ is finitely generated as a $\mathcal{B}$-module, and hence as an $\mathcal{A}$-module. 
Since $\overline{\mathcal{B}}$ is finitely generated as an $\mathcal{A}$-module, so is $C$.  

We claim that since $C$ is finitely generated as an $\mathcal{A}$-module, then $[C]_{s}=0$ for $s$ sufficiently large. 
Let $z_{1}, \ldots, z_{r}$ be a set of generators of $C$. 
Then since $H_{J}^{0}(C)=C$, there exists $k_{i}$ such that $J^{k_{i}}z_{i}=0$ for $1 \leq i \leq r$. Recall that $It \subset J$.
Therefore $[C]_{s}=0$ for $s > \max\{k_{i} + \deg(z_{i}) \mid 1 \leq i \leq r\}$.
Therefore, $\overline{I^{s}}+(x)/(x)= \overline{(I/(x))^{s}}$ for $s > \max\{k_{i} + \deg(z_{i}) \mid 1 \leq i \leq r\}$.
\end{proof}

Thus, by our main theorem (\cref{thm}) integral closure and specialization by a general element of an ideal $I$ are compatible for a large class of  ideals in many algebras over a field of characteristic zero. 
Moreover if the ring is normal then by \cref{powers}, integral closure of large enough powers of the ideal, $\overline{I^s}$ with $s >> 0$, is compatible with specialization by a general element of $I$. 
Based on computations, we believe that \cref{powers} may be true without requiring a normal ring, but different techniques would need to be used to prove such a result.

\section{Specialization by a General Linear Form}

In this section, we now shift our attention to the special case of monomial ideals in polynomial rings.
We consider the case in which we specialize the integral closure of an ideal by a general element of the unique homogeneous maximal ideal rather than a general element of the ideal $I$. 
The integral closure of an ideal does not behave as well with respect to specialization by a general element of the maximal ideal as it does with respect to a general element of the ideal. 
We give an example of an integrally closed monomial ideal of height $2$ which does not specialize with respect to a general linear form, even though by \cref{thm} it specializes with respect to a general element of the ideal itself.  

\begin{example}
Let $R = \mathbb{Q}[x,y,z]$. Let $I = \left(x^2, y\,z \right)$. Note that $I$ is an integrally closed height $2$ ideal of $R$. Let
\begin{equation*}
a = \alpha\,x+\beta\,y+\gamma\,z
\end{equation*}
with $\alpha, \beta, \gamma$ nonzero. 
We see below that $I+(a)/(a)$ is not an integrally closed ideal of $R/(a)$.

Notice that $z^{2}+(a)$ satisfies an equation of integral dependence over $I+(a)/(a)$ in $R/(a)$. 
Since $yz \in I$, $\beta\, y\,z \in I$. 
Therefore, $\alpha\, x\,z+\gamma\, z^{2} \in I+(a)$. 
One can verify that
\begin{align*}
    (z^{2})^{2}+\left(\frac{2\,\beta}{\gamma}y\,z-\frac{\alpha^{2}}{\gamma^{2}}\,x^{2}\right)(z^{2})+\frac{1}{\gamma^{2}}(\alpha x\,z + \gamma \,z^{2})^{2}  \in (a).
\end{align*}
Since $2\frac{\beta}{\gamma}yz-\frac{\alpha^{2}}{\gamma^{2}}x^{2} \in I \subseteq I+(a)$ and $\frac{1}{\gamma^{2}}(\alpha\,x\,z+\gamma\,z^{2})^{2}\in (I+(a))^{2}$, this shows that $z^{2}+(a) \in \overline{I+(a)/(a)}$.

One can verify that $z^{2}+(a) \not\in I+(a)/(a)$. 
Therefore, $z^{2}+(a) \in \overline{I+(a)/(a)} \setminus I+(a)/(a)$. 

\end{example}

More generally, there are classes of ideals whose integral closures do not specialize with respect to a general linear form. 
We give one such example in \cref{non-specialization-class}.

\begin{prop}\label{non-specialization-class}
Let $R=k[x_{1},\ldots,x_{d}]$ be a polynomial ring over an infinite field $k$.
Let $\mathfrak{m}=(x_{1},\ldots,x_{d})$ denote the homogeneous maximal ideal of $R$.
Let $I$ be an integrally closed ideal of height $d-1$ generated by $t$ forms of degree $n$ with $t<\binom{n+d-2}{n}$. 
Let $a=\alpha_{1}x_{1}+\cdots+\alpha_{d}x_{d}$ be a general linear form. 
Then $I+(a)/(a)$ is not an integrally closed $R/(a)$-ideal. 
\end{prop}

\begin{proof}

Since $a$ is general, we may assume $a \not\in I$. 

Since $\hgt{I}=d-1$, we have that 
\begin{equation*}
    \dim R/I=\dim R- \hgt{I}= 1.
\end{equation*} 
Hence $\dim R/(I+(a))=0$. Therefore $I+(a)/(a)$ is an $\mathfrak{m}/(a)$-primary ideal in $R/(a)$. 

Note that we have an isomorphism
\begin{equation*}
   \varphi: R/(a) \longrightarrow k[x_{1},\ldots,x_{d-1}],
\end{equation*}
defined by mapping 
\[x_{d} \mapsto \frac{1}{\alpha_{d}}(-\alpha_{1}x_{1}-\ldots-\alpha_{d-1}x_{d-1}).\]
Let $\mathfrak{n}=\left(x_{1},\ldots,x_{d-1}\right)$ denote the homogeneous maximal ideal of $k[x_{1},\ldots, x_{d-1}]$. Let $J=\varphi(I+(a)/(a))$. 
Since $I+(a)/(a)$ is $\mathfrak{m}/(a)$-primary, $J$ is $\mathfrak{n}$-primary. 
Moreover, $J$ is generated by forms of degree $n$. 

Suppose toward contradiction that $I + (a) / (a)$ is integrally closed. 
Then $J = \varphi(I + (a) / (a))$ is integrally closed. 
Since $J$ is generated by forms of degree $n$, $J \subseteq \mathfrak{n}^{n}$. 
Since $J$ is $\mathfrak{n}$-primary, there is $\mathfrak{n}^{s} \subseteq J \subseteq \mathfrak{n}^{n}$ for some $s \geq n$.
Therefore, there exists $t \geq 1$ such that
\begin{equation*}
    \mathfrak{n}^{nt} \subseteq \mathfrak{n}^{s} \subseteq J \subseteq \mathfrak{n}^{n}.
\end{equation*}
Since $J\subseteq \mathfrak{n}^{n}$, $J(\mathfrak{n}^{n})^{t-1} \subseteq (\mathfrak{n}^{n})^{t}$. 
Since $(\mathfrak{n}^{n})^{t} \subseteq J$, and $J$ is generated in degree $n$, $(\mathfrak{n}^{n})^{t} \subseteq J(\mathfrak{n}^{n})^{t-1}$. 
Hence $J(\mathfrak{n}^{n})^{t-1}=\mathfrak{n}^{nt}$. 
In other words, $J$ is a reduction of $\mathfrak{n}^n$, hence $J=\mathfrak{n}^{n}$ because $J$ is integrally closed. 

Therefore, the minimal number of generators of $J = \mathfrak{n}^n$ must be $\binom{n+d-2}{n}$. 
Since $I$ is generated by $t < \binom{n+d-2}{n}$ forms, so is  $J = \varphi(I+(a)/(a))$. 
This is a contradiction, therefore $I+(a)/(a)$ is not integrally closed. 
\end{proof}

However, there are cases in which going modulo a general element of the maximal ideal does preserve the property of being integrally closed, such as when $R/I$ is reduced and $\depth{R/I} \geq 2$.

\begin{prop}\label{depth2}
Let $(R,m)$ be a local excellent algebra over an infinite field $k$. 
Assume the residue field of $R$ is separable over $k$. Let $\mathfrak{m}=\left(x_{1},\ldots,x_{n}\right)$. 
Let $I$ be an $R$-ideal such that $R/I$ is reduced and $\depth{R/I} \geq 2$. 
Let $a= \sum_{i=1}^{n}\alpha_{i}x_{i}$ be a general element of $\mathfrak{m}$. 
Then $I+(a)/(a)$ is an integrally closed ideal of $R/(a)$. 
\end{prop}

\begin{proof}
We claim that $R/(I+(a))$ is reduced. 
Since $R/I$ is reduced and $a$ is a general element of the maximal ideal $\mathfrak{m}$, by \cite[Corollary~4.3]{Flenner}, $R/(I+(a))$ is reduced locally on the punctured spectrum $\Spec{R} \setminus \{\mathfrak{m}\}$. 

Since $\depth{R/I} \geq 2$ and $a$ is a general element of $\mathfrak{m}$, $\depth{R/(I+(a))} \geq 1$. 
Therefore, the maximal ideal $\mathfrak{m}/(I+(a))$ contains a nonzerodivisor on $R/(I+(a))$ and hence $\mathfrak{m}/(I+(a))$ is not an associated prime of $R/(I+(a))$. 
Since a ring is reduced if and only if it is reduced locally at associated primes, $R/(I+(a))$ is reduced. 
Therefore, $I+(a)/(a)= \sqrt{I+(a)/(a)}$, and hence
\[ I+(a)/(a) = \overline{I+(a)/(a)}.\]
\end{proof}

Next, we consider the case of squarefree monomial ideals in polynomial rings over an infinite field. 
Recall that a squarefree monomial ideal $I$ has a primary decomposition of the form 
\begin{equation*}
    I=\bigcap_{j=1}^{n}p_{j},
\end{equation*}
where $p_{j}=\left(x_{j_{1}},\ldots,x_{j_{k}}\right)$ with $1 \leq j_{1} < \ldots < j_{k} \leq d$ for each $1 \leq j \leq n$ by 
\cite[Corollary~1.3.4]{Herzog-Hibi}.
Therefore, every squarefree monomial ideal is an intersection of prime ideals, hence it is integrally closed. 

The following result shows that the integral closures of squarefree monomial ideals which are intersections of prime ideals generated by pairwise disjoint sets of variables specialize with respect to a general linear form. 

\begin{prop}\label{monomial}
Let $R=k[x_{1},\ldots,x_{d}]$ be a polynomial ring over an infinite field $k$.
Let $I=\bigcap_{i=1}^{n}p_{i}$, where each $p_{i}$ is generated by variables and the variable generating sets of $p_{i}$ and $p_{j}$ are disjoint for $i \not= j$. 
Let $a= \sum_{i=1}^{d}\alpha_{i}x_{i}$ be a general linear form. 
Then $I+(a)/(a)$ is an integrally closed ideal of $R/(a)$. 
\end{prop}

\begin{proof}
Since $I=\bigcap_{i=1}^{n}p_{i}$ and the ideals $p_{i}$ for $1 \leq i \leq n$ are generated by disjoint sets of variables, every monomial generator of $I$ is in $\prod_{i=1}^{n}p_{i}$. 
Since $\prod_{i=1}^{n}p_{i} \subseteq \bigcap_{i=1}^{n}p_{i}=I$, $I=\prod_{i=1}^{n}p_{i}$. 

We prove the result by induction on $n$. 

Base case: We first prove the result for $n=1$. 
Notice that if $I=p$, a prime ideal generated by variables, then $I+(a)$ is a prime ideal of $R$, and hence $I+(a)/(a)$ is a prime ideal of $R/(a)$. 
Since any prime ideal is integrally closed, $\overline{I+(a)/(a)}=I+(a)/(a)$.

Now let $n \geq 2$, and suppose the result holds for any intersection of $n-1$ primes generated by sets of variables which are pairwise disjoint.

Let $I=\bigcap_{i=1}^{n}p_{i} \subseteq \mathfrak{m}$, where $\mathfrak{m}$ is the maximal ideal of $R$. 

We consider two cases. 

Case 1: Suppose $\sum_{i=1}^{n}p_{i} \subsetneq \mathfrak{m}$. 
Since $I$ is equal to the intersection of its minimal primes, $R/I$ is reduced.

Next, we show that $\depth{R/I} \geq 2$.
Notice that $\dim{R/I} = \max\{\dim{R/p_{i}}\mid i=1,\ldots,n\}$.
Since $\sum_{i=1}^{n} p_i \subsetneq \mathfrak{m}$ and each $p_i \neq p_j$ are generated by nonempty disjoint sets of variables we have that $p_i \subsetneq p_i+(x_{j_{1}})) \subsetneq \mathfrak{m}$, where $x_{j_{1}}$ denotes a variable which is a generator of $p_{j}$ but not of $p_{i}$. 
These chains of primes ideals indicate that $\dim{R/p_i} \geq 2$ for $1 \leq i \leq n$, and hence $\dim{R/I} \geq 2$.

Without loss of generality, assume $x_{d} \not\in \sum_{i=1}^{n}p_{i}$. 
Since all generators in the squarefree monomial generating set of $I$ are not multiples of $x_{d}$,
the image of $I$ under the natural map $R \rightarrow R/(x_{d})$ is a squarefree monomial ideal with the same generators as $I$, and is therefore reduced. 
Since $\dim{R/I} \geq 2$ and hence $\dim{R/(I,x_{d})} \geq 1$, we conclude that $\depth{R/(I,x_{d})} \geq 1$. 
Therefore, $\depth{R/I} \geq 2$. 

Applying \cref{depth2}, we conclude that $I+(a)/(a)$ is integrally closed.

Case 2: Suppose $\sum_{i=1}^{n}p_{i} = \mathfrak{m}$. 
Let $\omega+(a) \in \overline{I+(a)/(a)}$. Since $I=\prod_{i=1}^{n}p_{i} \subseteq \prod_{i=1}^{n-1}p_{i}$, one has that 
\begin{equation*}
    \omega+(a) \in \overline{\prod_{i=1}^{n-1}p_{i}+(a)/(a)}.
\end{equation*}
By the inductive hypothesis,
\begin{equation*}
    \omega+(a) \in \prod_{i=1}^{n-1}p_{i}+(a)/(a).
\end{equation*}
Since $\omega+(a) \in \overline{I+(a)/(a)}$, there exists 
$m \in \mathbb{N}$ and $a_{i} \in I^{i}$ such that 
\begin{equation*}
    \omega^{m}+a_{1}\omega^{m-1}+\ldots+a_{m} \in (a).
\end{equation*}

Since $\sum_{i=1}^{n}p_{i}=\mathfrak{m}$, the general linear form 
\begin{equation*}
    a=b+c,
\end{equation*}
where $b \in \sum_{i=1}^{n-1}p_{i}$ and $c \in p_{n}$. 

If $n = 2$, we may assume that $\omega$ is homogeneous of degree at least $2$ since $(a)$ is a homogeneous ideal and $I$ is generated in degree $2$.

Notice that since $I+(a)/(a) \subseteq \overline{I+(a)/(a)}$, we may assume that $\omega= \omega_{1}+\omega_{2}$, in which $\supp{\omega_{1}} \subseteq p_1 \setminus p_2$ and $\supp{\omega_{2}} \subseteq p_2 \setminus p_1$.

We notice that
\begin{equation*}
    (\omega_{1}+\omega_{2})^{m}+a_{1}(\omega_{1}+\omega_{2})^{m-1}+\ldots+a_{m} \in (a),
\end{equation*}
modulo $p_2$ yields
\begin{equation*}
    (\omega_{1})^{m} \in (b) + p_2
\end{equation*}
Since $(b) + p_2$ is a prime ideal in $R$, then we see that 
\begin{equation*}
    \omega_{1} \in (b) + p_2.
\end{equation*}
Since $\supp{\omega_{1}} \subseteq p_1\setminus p_2$, we see that
\begin{equation*}
    \omega_{1} \in (b).
\end{equation*}
Moreover, by assumption $\omega_{1}$ has degree at least $2$ and therefore 
\begin{equation*}
    \omega_{1} \in (b) \cap p_1^{2} \subseteq (b)p_1.
\end{equation*}
Let $\omega_{1} = s b$ with $s \in p_1$. Then 
\begin{align*}
    \omega_{1}+(a) &=sb+ (a) \\
    &= -sc+(a) \\
    & \in I+(a).
\end{align*}
Repeating this argument modulo $p_1$ shows that $\omega_{2}+(a) \in I+(a)$ and therefore, $\overline{I+(a)/(a)}= I+(a)/(a)$.

If $n > 2$, without loss of generality, since $I+(a)/(a) \subseteq \overline{I+(a)/(a)}$, we may assume that 
\begin{equation*}
\supp{\omega} \subseteq \prod_{j=1}^{n-1}p_{j} \setminus p_{n}.
\end{equation*}

Then since 
\begin{equation*}
    \omega^{m}+a_{1}\omega^{m-1}+\ldots+a_{m} \in \left(a \right),
\end{equation*}
and each $a_{i} \in I \subseteq p_{n}$, we conclude that modulo $p_{n}$, 
\begin{equation*}
    \omega^{m} \in \left(b \right)+p_{n}.
\end{equation*}
Notice that $\left(b\right)+p_{n}$ is a prime ideal of $R$, and hence
\begin{equation*}
    \omega \in \left(b \right)+p_{n}.
\end{equation*}
Since $\supp{\omega} \subseteq \prod_{j=1}^{n}p_{j} \setminus p_{n}$, we conclude that 
\begin{equation*}
    \omega \in \left(b\right).
\end{equation*}
Let $s \in R$ be such that 
\begin{equation*}
    \omega= sb.
\end{equation*}
Then since we assume $\omega \in \prod_{j=1}^{n-1}p_{j}$ and $b \not\in p_{j}$ for any $1 \leq j \leq n-1$, we conclude that $s \in \prod_{j=1}^{n-1}p_{j}$. Therefore,
\begin{align*}
    \omega+(a) &= sb +(a) \\
    &= -sc +(a) \\
    &\in I+(a).
\end{align*}
This completes the proof that $\overline{I+(a)/(a)}=I+(a)/(a)$.
\end{proof}

We conjecture that the previous result holds more generally for squarefree monomial ideals without the condition of disjointness of generating sets of variables in the primary decomposition of the squarefree monomial ideal. 
Examples in the non-equigenerated case, such as $I=\left(x,y\right)\cap \left(x,z\right)=\left(x,yz\right)$ in $k[x,y,z]$, as well as in the equigenerated case, such as $I=\left(x,y\right)\cap \left(y,z\right) \cap \left(x,z\right)=\left(xy,xz,yz\right)$ in $k[x,y,z]$, support this conjecture.

\section{Acknowledgment}

The material for this paper comes from the authors' Ph.D. theses, under the direction of Professor Bernd Ulrich. We thank him for his suggestion of the problem, and for his invaluable support throughout the process.

\bibliographystyle{mandebib} 
\bibliography{references}

\end{document}